\documentclass[reqno]{amsart}
\usepackage{amssymb,latexsym}

\newtheorem{thm}{Theorem}[section]
\newtheorem{lemma}{Lemma}[section]
\newtheorem{cor}{Corollary}[section]
\newtheorem{prop}{Proposition}[section]

\theoremstyle{definition}
\newtheorem{defn}{Definition}[section]
\theoremstyle{remark}

\makeatletter
\@addtoreset{figure}{section}
\def\thefigure{\thesection.\@arabic\c@figure}
\def\fps@figure{h, t}
\@addtoreset{table}{section}
\def\thetable{\thesection.\@arabic\c@table}
\def\fps@table{h, t}
\@addtoreset{equation}{section}

\makeatother

\newcommand\hook{\mathbin{\hbox{\vrule height .5pt width 3.5pt depth 0pt
\vrule height 6pt width .5pt depth 0pt}}}

\begin{document}

\title[Multisymplectic Geometry, Covariant Hamiltonians, and Water Waves]
{Multisymplectic Geometry,  Covariant Hamiltonians,
and Water Waves}

\author[J.E. Marsden]{Jerrold E. Marsden}
\email{marsden@cds.caltech.edu}
\author[S. Shkoller]{Steve Shkoller}
\email{shkoller@cds.caltech.edu}
\address{Control and Dynamical Systems, 107-81,California Institute of Technology,
 Pasadena, CA 91125}
\date{June 25, 1997; To appear in Math. Proc. Camb. Phil. Soc.}

\begin{abstract}
This paper concerns the development and application of the multisymplectic
Lagrangian and Hamiltonian formalism for nonlinear partial differential
equations.
This theory generalizes and unifies the classical Hamiltonian formalism of
particle mechanics as well as the many pre-symplectic $2$-forms used by Bridges
[1997].  In this theory, solutions of a PDE are sections of a fiber bundle
$Y$ over
a base manifold $X$ of dimension $n$$+$$1$, typically taken to be
spacetime.  Given
a connection on $Y$, a covariant Hamiltonian density ${\mathcal H}$
is then intrinsically defined
on the primary constraint manifold $P_{\mathcal L}$, the image of the
multisymplectic version of the Legendre transformation. One views
$P_{\mathcal L}$
as  a subbundle of $J^1(Y)^\star$, the {\it affine} dual of $J^1(Y)$, the
first jet
bundle of $Y$.  A canonical multisymplectic ($n$$+$$2$)-form
$\Omega_{\mathcal H}$ is then defined, from which we obtain a
multisymplectic Hamiltonian system of differential equations that is
equivalent to both the original PDE as well as the Euler-Lagrange
equations of the corresponding Lagrangian.  Furthermore, we show that
the $n$$+$$1$ $2$-forms $\omega^{(\mu)}$ defined by Bridges [1997] are a
particular  coordinate representation for a single multisymplectic
($n$$+$$2$)-form, and in the presence of symmetries, can be {\it assembled}
into
$\Omega_{\mathcal H}$.  A generalized Hamiltonian Noether theory is then
constructed which relates the action of the symmetry groups lifted to
$P_{\mathcal L}$ with the conservation laws of the system.  These conservation
laws are defined by our generalized Noether's theorem which recovers the
vanishing of the divergence of the vector of $n$$+$$1$ distinct momentum mappings
defined in Bridges [1997] and, when applied to water waves, recovers
Whitham's conservation of wave action.  In our view,
the multisymplectic structure provides the natural setting for
studying dispersive wave propagation problems, particularly the
instability of water waves, as discovered by Bridges.  After developing
the theory, we show its utility in the study of periodic pattern formation
and wave instability.
\end{abstract}

\maketitle

\section{Introduction}
\label{sec1}
The canonical symplectic structure on the cotangent bundle of a given
configuration space provides a natural correspondence between Hamiltonian
vector
fields that govern the evolution of conservative ordinary differential
equations and the Hamiltonian functions which describe them. The setting of
tangent and cotangent bundles also provides a natural setting for the
Lagrangian
description of dynamics and the Legendre transformation that connects the
Lagrangian and Hamiltonian points of view.

In either case, the use of vector fields for the description of the
dynamics is
natural, because for ordinary differential equations there is a single
distinguished variable, time.  On the contrary, in systems of partial
differential
equations, solutions depend on multiple variables, usually spatial as well as
temporal, so one can make the case that a single vector field is not the
appropriate point of view because it would require collapsing all of the
spatial
structure of a solution to a  single point of phase space.  This occurs when a
choice is made to consider the time coordinate separately, and  describe the
dynamics in terms of an infinite-dimensional  space of fields at a given
instant
in time.  Although this methodology has been very successful, availing
itself to
the powerful organizing structure of the theory of evolution operators from a
point of view of functional analysis,  its immediate affect is a break of
manifest
covariance.

To maintain a covariant description, one can use a generalization of symplectic
geometry known as multisymplectic geometry.  This subject has a long
and distinguished history that we shall not review in this article;
rather, we follow the framework established in Gotay [1991], Gotay
et al [1992], and Gotay and Marsden [1992], wherein relativistic field theories
with Dirac-Bergmann type constraints are considered in a Lagrangian formalism,
while the Hamiltonian formalism relies on a ``space $+$ time'' (or
$3$$+$$1$) split. These references contain citations to much of the important
literature and history of the subject.

It is interesting that the structure of connection is
not necessary to intrinsically define the Lagrangian formalism (as shown in the
preceding references), while for the  intrinsic definition of a covariant
Hamiltonian the introduction of such a structure is essential. Of course, one
can avoid a connection if one is willing to confine ones attention to
local coordinates. We give an intrinsic definition of the covariant Hamiltonian
so that we may examine the fundamental interplay of the equivariance of
Hamiltonian and Lagrangian densities with respect to group actions.

Our objective is to use a variant of the multisymplectic Hamiltonian formalism
to generalize and make intrinsic the seminal and extremely important
work of Bridges [1996a, 1996b, 1997] on wave propagation, periodic
pattern formation, and linear instability.  Roughly speaking, our main result
states that in the case of $n$ distinct and possibly unbounded
spatial directions, the $n$$+$$1$ pre-symplectic $2$-forms introduced by
Bridges are actually contained in a single higher degree multisymplectic
$(n+2)$-form, and that in the presence of symmetries, these many forms can
be {\it assembled} into this single canonical form.  Furthermore, the
covariant Hamiltonian Noether theory that we construct, generalizes Bridges'
clever decomposition of water wave conservation laws, and
is an intrinsic restatement of the constrained variational
principles which lead to the existence of water wave instabilities
and diagonal periodic pattern formation.

We begin in Section \ref{sec2} by recalling some of the basic
constructions and a few key results from the multisymplectic formalism
of Gotay et al [1992].  In Section \ref{sec3}, we add the structure
of connection and intrinsically define our covariant Hamiltonian
density.  In Section \ref{sec4}, we show that our multisymplectic
formalism generalizes the classical theory of particle mechanics,
as well Bridges' theory of nonlinear partial differential equations.
Section \ref{sec5} is devoted to our development of a covariant
Hamiltonian Noether theory.  In Section \ref{sec6}, we show how this
theory recovers the classical conservation laws of particle
mechanics as well as the new conservation laws proposed by
Bridges [1996b] for studying water waves.  Finally, in Section
\ref{sec7} we show how our general theory applies to the study
of periodic pattern formation and the instability of waves.

\section{Multisymplectic Geometry}
\label{sec2}
A {\bf covariant configuration bundle} is a finite-dimensional fiber
bundle $\pi_{XY}:Y \rightarrow X$ over an oriented manifold $X$. In
many examples, especially those occurring in relativistic field
theories, $X$ is chosen to be spacetime and the fields of interest are
sections of this bundle. For nonrelativistic theories, such as
nonlinear waves, one typically chooses $X$ to be classical spacetime
(i.e., the product of the reals, ${\mathbb R}$, with the spatial
variables).

We shall need a little notation. Denote the fiber $\pi_{XY}^{-1}(x)$ of
$Y$ over $x\in X$ by $Y_x$ and the tangent space to $X$ at $x$ by
$T_xX$, etc., and denote sections of $\pi_{XY}$ by $\Gamma(\pi_{XY})$.
We also let $VY \subset TY$ be the {\bf vertical subbundle}; this is
the bundle  over $Y$ whose fibers are given by
\begin{equation}
V_yY = \{ v \in T_yY \mid T\pi_{XY} \cdot v  = 0 \},
\label{A2}
\end{equation}
where $T\pi_{XY} \cdot v$ denotes the derivative of the map $\pi_{XY}$
in the direction $v$.

Just as the covariant configuration bundle is the analogue of the
configuration space in particle mechanics, the first jet bundle,
defined next, is the field theoretic analogue of the tangent bundle.
\begin{defn}
The {\bf first jet bundle} $J^1(Y)$ is the {\it affine}
bundle over $Y$ whose fiber over $y\in Y_x$ consists of those linear
mappings $\gamma:T_xX \rightarrow T_yY$ satisfying
\begin{equation}
T\pi_{XY} \circ \gamma = \mbox{Identity on } T_xX. \qquad
\Diamond \label{A1}
\end{equation}
\end{defn}
The vector bundle underlying this affine bundle is the bundle whose
fiber over $y\in Y_x$ is the space $L(T_xX,V_yY)$ of linear mappings of
$T_xX$ into $V_yY$.  Note that for each $\gamma\in J^1(Y)_y$, we have
the  splitting
\begin{equation}
T_yY= \mbox{image }\gamma \oplus V_yY. \label{A3}
\end{equation}

The choice of the {\it first} jet bundle $J^1(Y)$ is used for the field
theoretic tangent bundle for classical field theories whose Lagrangians
depend on the point values of the fields and their {\it first}
derivatives.  For higher order field theories, one uses higher order
jet bundles; see Gotay et al [1992] and Gotay and Marsden [1992] for
references to this literature.

We let dim $X = n$$+$$1$ and the fiber dimension of $Y$ be $N$.
Coordinates on $X$ are denoted $x^\mu, \mu = 1, 2, \dots , n, 0$, and
fiber coordinates on $Y$ are denoted by $y^A, A = 1, \dots , N$.
These induce coordinates $v^A{}_\mu$ on the fibers of $J^1(Y)$.  If
$\phi : X \rightarrow Y$ is a section of $\pi_{XY}$, its tangent map at $x \in
X$, denoted $T_x\phi$, is an element of $J^1(Y)_{\phi(x)}$. Thus, the
map $x \mapsto T_x\phi$ defines a section of $J^1(Y)$ regarded as a bundle
over $X$.  This section is denoted $j^1(\phi)$ and is called the
{\bf first jet} of $\phi$.  In coordinates, $j^1(\phi)$ is given by
\begin{equation}
x^\mu \mapsto (x^\mu, \phi^A(x^\mu), \partial_\nu \phi^A (x^\mu)),
\label{A4}
\end{equation}
where $\partial_\nu = \partial / \partial x^\nu$.  A section of
the bundle $J^1(Y) \rightarrow X$ which is the first jet of a section of $Y
\rightarrow X$ is said to be {\bf holonomic}.

The field theoretic analogue of the cotangent bundle is defined next.
\begin{defn}
The {\bf dual jet bundle} $J^1(Y)^\star$ is the {\it vector} bundle
over $Y$ whose fiber at $y\in Y_x$ is the set of {\it affine} maps from
$J^1(Y)_y$ to $\Lambda^{n+1}(X)_x$, the bundle of $(n+1)$-forms on $X$.
\qquad $\Diamond$
\end{defn}
A smooth section of $J^1(Y)^\star$ is therefore an affine bundle map of
$J^1(Y)$ to $\Lambda^{n+1}(X)$ covering $\pi_{XY}$.  We choose affine maps
since $J^1(Y)$ is an affine bundle, and we map into $\Lambda^{n+1}(X)$ since
we are ultimately thinking of integration as providing the pairing on
sections.

Fiber coordinates on $J^1(Y)^\star$ are $(p, p_A{}^\mu)$, which
correspond to the affine map given in coordinates by
\begin{equation}
v^A{}_\mu \mapsto (p + p_A{}^\mu v^A{}_\mu) d^{n+1} x, \label{B1}
\end{equation}
where $d^{n+1}x = dx^1 \wedge \dots \wedge dx^n \wedge dx^0$.

Analogous to the canonical one- and two-forms on a cotangent bundle,
there are canonical forms on $J^1(Y)^\star$.  To define these, another
description of $J^1(Y)^\star$ will be convenient.  Namely,
let $\Lambda := \Lambda^{n+1}(Y)$ denote the bundle
of $(n+1)$-forms on $Y$, with fiber over $y \in Y$ denoted by
$\Lambda_y$ and with projection $\pi_{Y\Lambda}: \Lambda \rightarrow
Y$.  Let $Z \subset \Lambda$ be the subbundle whose fiber is given by
\begin{equation}
Z_y  =  \{z \in \Lambda_y \,\vert\, {v \hook ({w \hook z})} = 0 \text{
for all } v, w \in V_yY\},\label{B2}
\end{equation}
where ${v \hook \,\cdot \,}$ denotes left interior multiplication by $v$.

Elements of $Z$ can be be written uniquely as
\begin{equation}
z  =  pd^{n+1}x + p_A{}^\mu dy^A \wedge d^nx_\mu,
\label{B3}
\end{equation}
where
\[
d^n x_\mu  =  { \partial_\mu \hook d^{n+1}x } \quad\text{ and,
as before, }\quad \partial_\mu = \frac{\partial}{\partial x^\mu}.
\]
Hence, fiber coordinates for $Z$ are also $(p, p_A{}^\mu)$.

Corresponding to equating the coordinates $(x^\mu, y^A, p, p_A{}^\mu)$ of
$Z$ and
of $J^1(Y)^\star$, there is a vector bundle isomorphism
\begin{equation}
\Phi : Z \rightarrow J^1(Y)^\star. \label{B5}
\end{equation}
Intrinsically, $\Phi$ is defined by pull-back:
\begin{equation}
\Phi (z) (\gamma) = \gamma^*z \in \Lambda^{n+1} (X)_x \label{B6}
\end{equation}
where $z \in Z_y$, $\gamma \in J^1(Y)_y$ and $x = \pi_{XY}(y)$.  Using
fiber coordinates $v^A{}_\mu$ for $\gamma$,
the preceding equation becomes
\begin{equation}
\gamma^* dx^\mu  =  dx^\mu \quad \text{and} \quad \gamma^* dy^A  =  v^ A{}_
\mu
dx^\mu  \label{B7}
\end{equation}
and so
\begin{equation}
\gamma^*(pd^{n+1}x + p_A{}^\mu dy^A \wedge d^nx_\mu)  =
(p + p_A{}^\mu v^A{}_\mu)d^{n+1}x, \label{B8}
\end{equation}
where we have used $dx^\nu \wedge d^n x_\mu = \delta^\nu_\mu d^{n+1} x$.

One shows that the inverse of $\Phi$ can also be defined intrinsically,
although it is somewhat more complicated, and thus the spaces
$J^1(Y)^\star$ and $Z$ are canonically isomorphic as vector bundles
over $Y$.

There are canonical forms on $Z$ and the isomorphism
between $J^1(Y)^\star$ and $Z$ can be used to transfer these to
$J^1(Y)^\star$.  We first define the {\bf canonical} $(n+1)${\bf-form}
$\Theta_\Lambda$ on $\Lambda$ by
\begin{equation}
\aligned
 \Theta_\Lambda(z)(u_1, \dots, u_{n+1}) &= z(T\pi_{Y\Lambda} \cdot u_1, \dots,
T\pi_{Y\Lambda}
\cdot u_{n+1})\\ &=(\pi^*_{Y\Lambda}z)(u_1, \dots, u_{n+1})
\endaligned
 \label{B9}
\end{equation}
where $z \in \Lambda$ and $u_1, \dots, u_{n+1} \in T_z\Lambda$.  Define the
{\bf canonical} $(n+2)${\bf-form} $\Omega_\Lambda$ on $\Lambda$ by
\begin{equation}
\Omega_\Lambda = -d \Theta_\Lambda. \label{B10}
\end{equation}
Note that if $n = 0$ (i.e., $X$ is one-dimensional), then $\Lambda = T^*Y$
and $\Theta_\Lambda$ is the standard canonical one-form.  If $i_{\Lambda Z}: Z
\rightarrow
\Lambda$  denotes the inclusion, the {\bf canonical} $(n+1)${\bf-form}
$\Theta$ on $Z$ is defined by
\begin{equation}
\Theta = i^*_{\Lambda Z} \Theta_\Lambda \label{B11}
\end{equation}
and the {\bf canonical} $(n+2)${\bf-form} $\Omega$ on $Z$ is defined by
\begin{equation}
\Omega = -d \Theta = i^*_{\Lambda Z} \Omega_\Lambda. \label{B12}
\end{equation}
The pair $(Z, \Omega)$ is called {\bf multiphase space} or
{\bf covariant phase space}. It is an example of a {\bf multisymplectic
manifold}.

Using (\ref{B3}), (\ref{B9}), (\ref{B10}), (\ref{B11}) and (\ref{B12}),
one finds that the coordinate expression  for $\Theta$ is
\begin{equation}
\Theta = p_A{}^\mu dy^A \wedge d^nx_\mu + pd^{n+1} x, \label{B13}
\end{equation}
and so
\begin{equation}
\Omega = dy^A \wedge dp_A{}^\mu \wedge d^nx_\mu - dp \wedge d^{n+1} x.
\label{B14}
\end{equation}

Let the {\bf Lagrangian density} ${\mathcal L} : J^1(Y) \rightarrow
\Lambda^{n+1}(X)$, be a given smooth bundle map over $X$.  In coordinates, we
write
\begin{equation}
{\mathcal L} (\gamma)  =  L(x^\mu, y^A, v^A{}_\mu) d^{n+1}x. \label{C1}
\end{equation}

The corresponding {\bf covariant Legendre transformation} associated with
${\mathcal L}$ is a fiber preserving map over $Y$, ${\mathbb F}{\mathcal L}:
J^1(Y) \rightarrow J^1(Y)^\star \cong  Z $, whose intrinsic definition
follows.
\begin{defn}
If $\gamma \in J^1(Y)_y$, we define ${\mathbb F}{\mathcal L}(\gamma)\in
J^1(Y)^\star_y$ (where $y \in  Y_x$) to be the first order vertical Taylor
approximation to ${\mathcal L}$:
\begin{equation}
{\mathbb F}{\mathcal L}(\gamma) \cdot  \gamma'  = {\mathcal L}(\gamma)
+\left.
\frac{d}{d\varepsilon} \right|_{\varepsilon = 0}{\mathcal L}(\gamma +
\varepsilon(\gamma' - \gamma)) \label{C3}
\end{equation}
where $\gamma' \in  J^1(Y)_y$. \qquad $\Diamond$
\end{defn}
A straightforward calculation shows that the covariant Legendre
transformation is given in coordinates by
\begin{equation}
p_A{}^\mu =
\frac{\partial L}{\partial v^A{}_\mu}, \text{ and }\quad p =
L -\frac{\partial L}{\partial v^A{}_\mu} v^A{}_\mu .\label{C2}
\end{equation}
Notice that formally, the second of these equations defines the
(negative of the) energy while the first one is reminiscent of the
usual relation $p_i = \partial L / \partial \dot{q}^i$  from classical
mechanics.
One of the nice features of the covariant Legendre transformation is how these
two basic aspects of the Legendre transformation arise from a single
construction.
\begin{defn} The {\bf Cartan form} is the $(n+1)$-form $\Theta_{{\mathcal
L}}$ on $J^1(Y)$  defined by
\begin{equation}
\Theta_{{\mathcal L}} =  (\mathbb F{\mathcal L})^*\Theta \label{D1}
\end{equation}
where $\Theta$ is the canonical $(n+1)$-form on $Z$.  We  also define
the $(n+2)$-form $\Omega _{\mathcal L}$ by
\begin{equation}
\Omega_{\mathcal L} = -d \Theta_{\mathcal L} = (\mathbb F {\mathcal L})^*
\Omega \label{D2}
\end{equation}
where $\Omega  = -d\Theta$ is the canonical $(n+2)$-form on $Z$.
$\Diamond$ \end{defn}

\section{The Covariant Hamiltonian}
\label{sec3}

In this section we develop an intrinsic covariant (or multisymplectic)
Hamiltonian formalism. We begin by noting that the covariant Legendre
transformation ${\mathbb F}{\mathcal L}:J^1(Y)\rightarrow J^1(Y)^\star$
is never a fiber bundle diffeomorphism since
dim$J^1(Y)^\star$$=$dim$J^1(Y)$$+$$1$; nevertheless, it is
appropriate in many examples to require ${\mathbb F}{\mathcal L}$ to be a
smooth bundle diffeomorphism over $Y$ onto its image.  In fact, from the
second equation in (\ref{C2}), the image of ${\mathbb F}{\mathcal L}$
defines the primary constraint of the theory.

\begin{defn}
We say that ${\mathcal L}$ is {\bf regular} if the image of the first
jet bundle under the covariant Legendre transformation
$P _{\mathcal L} := {\mathbb F}{\mathcal L}( J ^1 (Y)) $ is a smooth
manifold and if ${\mathbb F}{\mathcal L}$ is a diffeomorphism onto
$P _{\mathcal L}$. We call $P _{\mathcal L}$ the {\bf primary
constraint manifold}. \qquad
$\Diamond$
\end{defn}

One should note that many field theories, such as the vacuum Maxwell equations
and many others, especially relativistic ones, are not regular because of the
presence of constraints (such as $ {\rm div}\, {\bf E} = 0$ for Maxwell's
equations). We are assuming regularity only for simplicity and because it is
appropriate for the examples we have in mind. Gotay et al [1992] deal with the
more general case in a Lagrangian formalism.

At this point,  we introduce the additional structure of a connection.
While connections are not particularly needed for the Lagrangian side of
field theory, they seem to be essential for the development of an intrinsic
Hamiltonian formalism.   We recall the definition of an (Eheresmann)
connection
as a vertical-valued one-form.

\begin{defn}
A {\bf connection} on $Y$ is a vector bundle map ${\mathfrak A}:TY
\rightarrow VY$ such that on each fiber over
$y\in Y$,  ${\mathfrak A}:T_yY \rightarrow V_yY$ satisfies
\begin{equation}
{\mathfrak A} = \text{ Identity on } VY.  \label{A1'}
\end{equation}
The {\bf horizontal space} at each point $y \in Y $ is defined by
${\rm hor}_y = {\rm ker}\,{\mathfrak A}_y$, so that
we have $T _y Y = {\rm hor}_y \oplus V _y Y $. $\Diamond$
\end{defn}
In coordinates, the action of ${\mathfrak A}$ on a tangent vector to
$Y$, namely $(v ^\nu, v ^A)$ is written as $ (0, v ^A + {\mathfrak
A}_\mu^A  v ^\mu)$. This defines the coordinate expression for the
connection.
We remark that it is not entirely necessary to a priori explicitly introduce
a connection if one wishes to define the Hamiltonian locally in a coordinate
chart, and then use coordinate patches to obtain a global characterization;
however, the process of producing a coordinate independent global definition
is tantamount to producing a connection.

Next, we reexpress the covariant Legendre transformation
${\mathbb F}{\mathcal L}$ in  terms of a vertical derivative of
functions on $J^1(Y)$.

\begin{defn}
For any $x\in X$ and $y\in Y_x$,
let $U\subset J^1(Y)_y$ be an open subset, and let ${\mathcal S}\in
C^1(U,\Lambda^{n+1}(X)_x)$.  Then the {\bf covariant derivative} of $S$
associated with the connection ${\mathfrak A}$ maps $U$ into
$C^0(J^1(Y)_y,\Lambda^{n+1}(X)_x)$ and is defined, for any
$\gamma \in U$, and $\gamma'\in J^1(Y)_y$, by
\begin{equation}
D_{\mathfrak A}{\mathcal S}(\gamma) \cdot \gamma' =
\left.\frac{d}{d\varepsilon}
\right|_{\varepsilon=0}{\mathcal S}(\gamma+\varepsilon{\mathfrak A}
(\gamma')).  \qquad \Diamond \label{E00}
\end{equation}
\end{defn}
It is then natural to consider the covariant derivative of the smooth
bundle map ${\mathcal L}:\pi_{X,J^1(Y)}\rightarrow \pi_{X,\Lambda^{n+1}(X)}$,
so that
using (\ref{E00}), the Legendre transformation can be written as
\begin{equation}
{\mathbb F}{\mathcal L} (\gamma)\cdot \gamma' = \left[{\mathcal L} (\gamma)
-D_{\mathfrak A}{\mathcal L}(\gamma)\cdot \gamma\right] +
D_{\mathfrak A} {\mathcal L}(\gamma) \cdot \gamma' \label{E000}
\end{equation}
for all $\gamma, \gamma' \in J^1(Y)_y$ and $y\in Y_x$. Note that this
expression is affine with the first two terms being the constant
terms and the last one being the linear term.

We may now define the covariant Hamiltonian density on the primary
constraint manifold $P_{\mathcal L}$.
\begin{defn}
For a regular Lagrangian, the corresponding {\bf covariant Hamiltonian}
${\mathcal H}: P_{\mathcal L} \rightarrow \Lambda^{n+1}(X)$
is defined by
\begin{equation}
{\mathcal H} (z)= D_{\mathfrak A}{\mathcal L}(\gamma) \cdot \gamma
- {\mathcal L}(\gamma)
\label{D2b}
\end{equation}
where $z={\mathbb F}{\mathcal L}(\gamma)$.
$\Diamond$
\end{defn}
In coordinates, we may write ${\mathcal H} = H d^{n + 1}x$ where
\begin{equation}
H =  \frac{\partial L}{\partial v^A{}_\mu} (v^A{}_\mu +
{\mathfrak A}_\mu^A) - L.
\nonumber
\end{equation}
Notice that the covariant Hamiltonian is well defined under the
assumption of regularity; namely, the map $\gamma \mapsto z = {\mathbb
F}{\mathcal L}( \gamma ) $ from $J^1(Y)$ to $P _{\mathcal L}$ is a
diffeomorphism.

We coordinatize the primary constraint manifold by
$(x ^\mu, y ^A, p_A {}^\mu)$ with $p$ now
expressed in terms of the other variables by rewriting the preceding
expression for $H$ as
\begin{equation}
H =  \frac{\partial L}{\partial v^A{}_\mu}{\mathfrak A}_\mu^A  - p,
\nonumber
\end{equation}
regarded as an implicit equation for $p$.

Let $i_{J^1(Y)^\star,P_{\mathcal L}}:P_{\mathcal
L}\rightarrow J^1(Y)^\star$ denote the inclusion.  We may pull-back the
canonical $(n+1)$- and $(n+2)$-forms on $J^1(Y)^\star$ to $P_{\mathcal
L}$  and obtain (using a notation to remind us that this takes the
Hamiltonian point of view):
\begin{equation}
\begin{array}{c}
\Theta_{\mathcal H}= i_{J^1(Y)^\star,P_{\mathcal L}}^* \Theta, \\
\Omega_{\mathcal H}= i_{J^1(Y)^\star,P_{\mathcal L}}^* \Omega.
\end{array}
\label{N2}
\end{equation}
In canonical coordinates, we have
\begin{equation}
\begin{array}{c}
\Theta_{\mathcal H} = {p_A}^\mu dy^A\wedge d^nx_\mu +({p_A}^\mu
{\mathfrak A}^A_\mu - H)\wedge d^{n+1}x
\\
\Omega_{\mathcal H} = dy^A \wedge {dp_A}^\mu \wedge d^nx_\mu +
\left[\frac{\partial H}{\partial y^A}dy^A
+ \left(\frac{\partial H}{\partial {p_A}^\mu}-{\mathfrak A}^A_\mu
\right) d{p_A}^\mu\right] \wedge d^{n+1}x.
\end{array}
\label{E0}
\end{equation}

For many important examples, we will consider $X$ as the classical
spacetime manifold with the locally trivial connection which is simply
the natural projection whose action in coordinates is $(0,v^A)$, i.e.,
the components ${\mathfrak A}^A_\mu =0$.

\begin{defn}
Let $\phi \in \Gamma(\pi_{XY})$, and $j^1(\phi)$ its first jet.
A section $z$ of $\pi_{X,P_{\mathcal L}}$  is called {\bf conjugate} to
$j^1(\phi)$
if $z = {\mathbb F}{\mathcal L} \circ j^1(\phi)$. In this case, we shall write
$\tilde{j}^1(\phi)$ for $z$ and say that $z$ is holonomic.
 $\Diamond$
\end{defn}

\begin{defn}
A holonomic section $z$ of $P_{\mathcal L}$ is called {\bf Hamiltonian}
for ${\mathcal H}$ if
\begin{equation}
z^* ({ U\hook  \Omega_{\mathcal H}}) = 0. \label{E1}
\end{equation}
for any $U\in T(P_{\mathcal L})$. We also refer to the system of equations
(\ref{E1}) regarded as differential equations for $z$ as the
{\bf multihamiltonian system} of equations associated to ${\mathcal H}$.
$\Diamond$
\end{defn}
\begin{lemma}
\label{LEMMA1}
If ${\mathbb F}{\mathcal L}:J^1(Y)\rightarrow P_{\mathcal L}$ is a fiber
bundle diffeomorphism
over $Y$ and $\phi \in \Gamma(\pi_{XY})$, then the following are equivalent:
\begin{itemize}
\item[{\bf (i)}]$\tilde{j}^1(\phi)^* (U \hook\Omega_{\mathcal H}) =0$ for
any $U\in
T(P_{\mathcal L})$;
\item[{\bf (ii)}]$j^1(\phi)^* (W \hook \Omega_{\mathcal L}) =0$ for any $W\in
T(J^1(Y))$.
\end{itemize}
\end{lemma}
\begin{proof}
Assume {\bf (i)} holds and let $U\in T(P_{\mathcal L})$.  Since ${\mathbb
F}{\mathcal L}$ is a fiber bundle diffeomorphism, there exists $W \in
T(J^1(Y))$
such that $T{\mathbb F}{\mathcal L}
\circ W = U \circ {\mathbb F}{\mathcal L}$.  Hence,
\begin{eqnarray*}
0=\tilde{j}^1(\phi)^* (U\hook\Omega_{\mathcal H})& =& j^1(\phi)^* {\mathbb
F}{\mathcal L}^* ({T{\mathbb F}{\mathcal L} \cdot W}\hook\Omega_{\mathcal
H})\\ &=&
j^1(\phi)^* ({{\mathbb F}{\mathcal L}^*T{\mathbb F}{\mathcal L}
\cdot W}\hook{\mathbb F}{\mathcal L}^*\Omega_{\mathcal H})\\ &=& j^1(\phi)^*
(W\hook\Omega_{\mathcal L}).
\end{eqnarray*}
Using the same argument, the inverse function theorem guarantees that
the converse holds as well.
\end{proof}
We are thus led to the following conclusion.
\begin{thm}
\label{THM1}
If ${\mathbb F}{\mathcal L}:J^1(Y)\rightarrow P_{\mathcal L}$ is a fiber
bundle diffeomorphism
over $Y$ and $\phi \in \Gamma(\pi_{XY})$, then the following are equivalent:
\begin{itemize}
\item[{\bf (i)}]$\phi$ is a stationary point of $\int_X{\mathcal L}
(j^1(\phi))$;
\item[{\bf (ii)}]$\tilde{j}^1(\phi)$ is a Hamiltonian section for
${\mathcal H}$.
\end{itemize}
\end{thm}
Before we prove the theorem, we state the following definition and
lemma.
\begin{defn} A {\bf (finite) variation} of $\phi$ is a curve
$\phi_\lambda  = \eta_\lambda \circ \phi$,  where $\eta_\lambda $ is the flow
of a  vertical vector field $V$ on $Y$ which is compactly supported in $X$.
One says that $\phi$ is a {\bf stationary point} of the action if
\begin{equation}
\left.\frac{d}{d\lambda} \left[ \int_X {\mathcal L}(j^1(\phi_\lambda)) \right]
\right|_{\lambda=0} = 0 \label{E7}
\end{equation}
for all variations $\phi_\lambda$ of $\phi$.
\end{defn}
\begin{lemma}
\label{LEMMA2}
If ${\mathbb F}{\mathcal L}:J^1(Y)\rightarrow P_{\mathcal L}$ is a fiber
bundle diffeomorphism, then
\[
\tilde{j}^1(\phi)^*( { U \hook \Omega_{\mathcal H}}) =0
\]
for any $U$ which is $\pi_{Y,P_{\mathcal L}}$-vertical or is tangent to
$\tilde{j}^1(\phi)$. Similarly,
\[ j^1(\phi)^*( { W \hook \Omega_{\mathcal L}}) =0
\]
for any $W$ which is $\pi_{Y,J^1(Y)}$-vertical or is tangent to $j^1(\phi)$.
\end{lemma}
\begin{proof}
Since ${\mathbb F}{\mathcal L}$ is a {\it fiber-preserving} bundle
diffeomorphism,
for any
$\pi_{Y, P_{\mathcal L}}$-vertical $U$, there exists a
$\pi_{Y,J^1(Y)}$-vertical
$W$ such that $T{\mathbb F}{\mathcal L} \circ W = U \circ {\mathbb
F}{\mathcal L}$.
Using canonical coordinates, let us write $U$ and $W$ as
\begin{equation}
U = {U_A}^\mu \frac{\partial}{\partial {p_A}^\mu }, \mbox{ and }
W={W_\mu}^A \frac{\partial}{\partial {v^A}_\mu}.
\nonumber
\end{equation}
A calculation using (\ref{E0}) shows that
\begin{equation}
U \hook \Omega_{\mathcal H} =  {U_A}^\mu \left( dy^A\wedge d^n x_\mu +
\left(\frac{ \partial H} {\partial {p_A}^\mu} -{\mathfrak A}^A_\mu
\right) d^{n+1}x\right).
\nonumber
\end{equation}
Hence,  using Lemma \ref{LEMMA1}, we have that
\begin{equation}
\tilde{j}^1(\phi)^* ({U \hook \Omega_{\mathcal H}})  =
j^1(\phi)^* \left(-W_{{v_\nu}^B} \frac{\partial L}{\partial {v^A}_\mu
\partial {v^B}_\nu} (dy^A\wedge d^nx_\mu - {v^A}_\mu d^{n+1}x )\right),
\end{equation}
which vanishes using (\ref{B7}).
On the other hand, if $U$ is tangent to the graph of $\tilde{j}^1(\phi)$, then
$U = T\tilde{j}^1(\phi) \cdot v$ for some $v \in TX$ so that
\begin{equation}
\tilde{j}^1(\phi)^* ({U \hook \Omega_{\mathcal H}})  =
\tilde{j}^1(\phi)^*({(T\tilde{j}^1(\phi)
\cdot v) \hook \Omega_{{\mathcal H}}})  =  { v \hook
(\tilde{j}^1(\phi)^*\Omega_{{\mathcal H}})},
\nonumber
\end{equation}
which vanishes since $\tilde{j}^1(\phi)^*\Omega_{\mathcal L}$ is an $(n +
2)$-form on the
$(n + 1)$-manifold $X$. The identical argument works for $W$.
\end{proof}
The proof of Lemma \ref{LEMMA2} shows that in canonical coordinates
\begin{equation}
 \frac{\partial H} {\partial {p_A}^\mu} = {\nu^A}_\mu + {\mathfrak A}^A_\mu,
\label{E101}
\end{equation}
and in the case that $U = U_{y^A} \frac{\partial}{\partial y^A}$ and
that $\tilde{j}^1(\phi)^* ( { U \hook  \Omega_{\mathcal H}}) = 0$, we obtain
that
\begin{equation}
\frac{\partial H} {\partial y^A} = -\frac{\partial {p_A}^\mu}
{\partial x^\mu}. \label{E102}
\end{equation}
Thus, equations (\ref{E101}) and (\ref{E102}) are the coordinate
expressions for a multihamiltonian system.

\begin{proof}[Proof of Theorem \ref{THM1}]
Let $\phi_\lambda = \eta_\lambda \circ \phi$ of $\phi$  be a variation
corresponding to a $\pi_{XY}$-vertical vector field $V$ on $Y$ with
compact support in $X$. Using (\ref{B3}) we find that ${\mathcal
L}(j^1(\phi)) =
j^1(\phi)^*\Theta_{{\mathcal L}}$ and hence
\begin{equation}
\aligned
\left.\frac{d}{d\lambda} \left[ \int_X {\mathcal L}(j^1(\phi_\lambda)) \right]
\right|_{\lambda =0} & = \frac{d}{d\lambda} \int_Xj^1(\phi_\lambda)^*
\Theta_{\mathcal L} |_{\lambda=0} \\ & = \frac{d}{d\lambda}\left. \left[
\int_X j^1(\phi)^* j^1 (\eta_\lambda)^*
\Theta_{\mathcal L})\right] \right|_{\lambda=0} \\
%& = \int_Xj^1(\phi)^*\text{\bf{\it \$}}_{j^1(V)} \Theta_{\mathcal L}
& = \int_Xj^1(\phi)^*{\bf{\mathfrak L}}_{j^1(V)} \Theta_{\mathcal L}
\endaligned
\label{E9}
\end{equation}
where
\begin{equation}
j^1(V) =  \left.\frac{d}{d\lambda}  j^1(\eta_\lambda) \right|_{\lambda=0}
\nonumber
\end{equation}
is the jet prolongation of $V$ to $J^1(Y)$
(see Definition \ref{prolong}, if necessary).
Using Cartan's magic formula, we get ${\bf {\mathfrak
L}}{}_W\Theta_{\mathcal L}  =
-{ W \hook \Omega_{\mathcal L}} + d ( { W\hook \Theta_{\mathcal L}})$, which,
together with (\ref{E9}), gives
\begin{equation}
\aligned
\left.
\frac{d}{d\lambda}
\left[
\int_X {\mathcal L}(j^1(\phi_\lambda))
\right]\right|_ {\lambda=0} &
= -\int_Xj^1(\phi)^*({j^1(V) \hook
\Omega_{\mathcal L}}) \\ & \qquad  \qquad
+ \int_X  d j^1(\phi)^*({j^1(V) \hook
\Theta_{\mathcal L}}) \\
& = -\int_Xj^1(\phi)^*({j^1(V) \hook \Omega_{\mathcal L}})
\endaligned
\label{E10}
\end{equation}
by Stokes' theorem and the fact that $V$, and hence $j^1(V)$ is
compactly supported in $X$. Lemma \ref{LEMMA1} together with (\ref{E10})
shows that {\bf (ii)} implies {\bf (i)}.

The converse follows from the fact that any $\pi_{X,J^1(Y)}$-vertical
vector field $W$ may be decomposed as
\begin{equation}
W = j^1(V) + W_1, \nonumber
\end{equation}
where $V$ is $\pi_{XY}$-vertical and $W_1$ is $\pi_{Y,J^1(Y)}$-vertical.  Then
if {\bf (i)} holds, Lemma \ref{LEMMA1} and  Lemma \ref{LEMMA2} together
with (\ref{E10}) show that
\begin{equation}
\int_X \tilde{j}^1(\phi)^* ({ U \hook \Omega_{\mathcal H}}) =0 \nonumber
\end{equation}
for all vector fields $U$ on $P_{\mathcal L}$ with compact support in $X$.
Since the space of smooth vector fields on $J^1(Y)$ is a module over the
ring of smooth functions on $X$, an argument like that in the
Fundamental Lemma of the Calculus of Variations shows that the
integrand must vanish for all vector fields $U\in T(P_{\mathcal L})$ with 
compact support in $X$.  A partition of unity argument then shows that 
{\bf (ii)} implies {\bf (i)}.
\end{proof}

In the next section, we shall demonstrate the machinery of our intrinsic
development on two examples: classical mechanics and nonlinear PDEs.
We note that the essence of both of the following examples are the equations
(\ref{E101}) and (\ref{E102}).  For a development of a generalized
Hamiltonian structure based on Hamiltonian vector fields that seems
well-suited for ordinary differential equations, we refer the reader to
Cantrijn et al [1997].

\section{Particle Mechanics and Nonlinear PDEs}
\label{sec4}
In this section, we show that our multisymplectic Hamiltonian formalism
generalizes classical particle mechanics and is a natural setting for
nonlinear Hamiltonian partial differential equations.
\subsection{Particle mechanics}
For non-relativistic classical mechanics with a configuration manifold $Q$ (of
dimension $N$), we choose $X={\mathbb R}$ (so that $n = 0$) and $Y={\mathbb R}
\times Q$. In this case, $J^1(Y) = {\mathbb R} \times TQ$,  and the
cross-product induces a (flat) connection  ${\mathfrak A}  : {\mathbb R}
\times TQ \rightarrow TQ$.
The dual jet bundle is given by
$J^1(Y)^\star = T^*{\mathbb R}\times T^*Q$ and has canonical
coordinates $(t,p, q^1,...,q^N,p_1,...,p_N)$.

Given a Lagrangian in the usual sense
$L : TQ \rightarrow {\mathbb R}$, we define
${\mathcal L}: J ^1(Y) \rightarrow \Lambda ^1(Y)$ by
\[
{\mathcal L} (t,  q ^A, \dot{q} ^A ) = L (t, q ^A, \dot{q} ^A ) dt .
\]
The covariant Legendre transformation is the map ${\mathbb F}{\mathcal L}:
J ^1(Y)
\rightarrow J ^1 (Y) ^\ast$; i.e.,
${\mathbb F} {\mathcal L}: {\mathbb R} \times T Q \rightarrow T ^\ast
{\mathbb R}
\times T^\ast Q $ given by
\[
{\mathbb F} {\mathcal L} (t, q ^i, \dot{q} ^i) = ( t, L(t,q^A,\dot{q}^A)-
p _A \dot{q} ^A  , q ^A,  p _A)
\]
where $p _A = \partial L / \partial \dot{q}^A$.
In this case, it is clear that $P_{\mathcal L}= {\mathbb R} \times T^*Q$
(as a subbunble of $T^*{\mathbb R} \times T^*Q$) with coordinates
$(t,q^A,p_A)$.  Assume that the Legendre
transformation is nondegenerate in the usual sense so that
${\mathbb F} {\mathcal L} : {\mathbb R}\times TQ \rightarrow P_{\mathcal L}$
is a vector bundle diffeomorphism over ${\mathbb R}$ and the
corresponding Hamiltonian $H : T^\ast Q \rightarrow {\mathbb R} $ is well
defined.
The function $H$ corresponds to the density
${\mathcal H}: P_{\mathcal L}\rightarrow \Lambda ^1(X)$ as
\[
{\mathcal H} ( t, q ^A, p _A ) = H (t,q^A, p _A) dt ,
\]
where
\[
{\mathcal H}(z)  =   D_{\mathfrak A}{\mathcal L} (\gamma ) \cdot \gamma
- {\mathcal L}( \gamma ) , \ \
z = {\mathbb F}{\mathcal L}( \gamma)\in P_{\mathcal L}
\]
has  the coordinate expression
\[
{\mathcal H}( t, q ^A, p _A ) = (p_A q^A- L(t,q^A,\dot{q}^A))dt.
\]

In this case, we obtain the usual symplectic $2$-form on extended phase space
\begin{equation}
\Omega_{\mathcal H} = dq^A\wedge dp_A +dH\wedge dt. \nonumber
\end{equation}
\begin{prop}
\label{PROP1}
Let $\bar{V}$ be a vector field on $P_{\mathcal L}$ with integral curve
$\tilde{j}^1(\phi)$. Then $\bar{V}$ is a Hamiltonian vector field for $H$
if and only if
$\tilde{j}^1(\phi)$ is a Hamiltonian section for ${\mathcal H}$.
\end{prop}
\begin{proof}

In coordinates, let $\bar{v}=(1,V_q^A,V_p^B)$ and let
$(U_t,U)\equiv(U_t,U_q^A,U_p^B)$ be the coordinates for an arbitrary
vector field $\bar{U}\in T(P_{\mathcal L})$.
Then
\begin{equation}
i_{\bar{U}}\Omega_{\mathcal H} = U_p^B dq^A - U_q^A dp^B +(dH \cdot
U)dt - U_t dH,
\label{E103}
\end{equation}
and the pull-back of (\ref{E103}) under $\tilde{j}^1(\phi)$ vanishes if and
only if
\begin{equation}
\dot{q}^A =-\frac{\partial H}{\partial p_A}, \qquad
\dot{p}_B =\frac{\partial H}{\partial q^B},
\end{equation}
and
\begin{equation}
dH \cdot  U  =0.
\end{equation}
\end{proof}
\subsection{Nonlinear partial differential equations}
\label{subsec_pde}
To motivate the exposition, consider the
nonlinear wave equation given by
\begin{equation}
\frac{\partial ^2 \phi}{\partial {x^0}^2} - \triangle \phi - V'(\phi)=0,
\ \ \phi \in \Gamma(\pi_{XY}),
\label{Z1}
\end{equation}
where $\triangle$ is the Laplace-Beltrami operator and $V$ is a real-valued
$C^\infty$ function of one variable.

We will show that equation (\ref{E1}) along the holonomic section
$\tilde{j^1}(\phi)$ is equivalent to the nonlinear PDE (\ref{Z1}) as well as
the  Bridges [1997] coordinate representation. For clarity of presentation, we
will consider only one spatial dimension. In terms of our general notation, we
set
$X={\mathbb R} ^2$ ($n=1$) and $Y={\mathbb R}^2 \times M$, so that
sections of $Y$ have the coordinate  expressions $(x^0,x^1,\phi)$,
and $TY = {\mathbb R}^2 \times TM$.  The cross-product once again induces
a flat connection
defined by the natural projection ${\mathfrak A}: {\mathbb R}^2 \times TM
\rightarrow TM$.

The first jet bundle $J^1(Y)$ is a five-dimensional manifold and
sections of $J^1(Y)$ have local coordinates $(x^0,x^1,\phi ,
\frac{\partial \phi}{\partial x^0}, \frac{\partial \phi}{\partial x^1})$.
The affine  dual $J^1(Y)^\star$ is six-dimensional with its sections
having the local coordinates $(x^0,x^1,\phi,p,p^0,p^1)$.
The Lagrangian density ${\mathcal L}:J^1(Y)\rightarrow \Lambda^2(X)$ is
expressed as
\[
{\mathcal L}\left( x^0,x^1,\phi ,
\frac{\partial \phi}{\partial x^0}, \frac{\partial \phi}{\partial x^1}\right)
= L\left( x^0,x^1,\phi
\frac{\partial \phi}{\partial x^0}, \frac{\partial \phi}{\partial x^1}\right)
dx^1 \wedge dx^0
\]
which, for the case of the nonlinear wave equation is
\[
{\mathcal L}\left( x^0,x^1,\phi ,
\frac{\partial \phi}{\partial x^0}, \frac{\partial \phi}{\partial x^1}\right)
= \left[\frac{1}{2}(\{\partial_0 \phi\}^2 - \{\partial_1 \phi\}^2)
+ V(\phi)\right] dx^1 \wedge dx^0. \]
In this setting, the covariant Legendre transformation ${\mathbb
F}{\mathcal L}:
J^1(Y)\rightarrow J^1(Y)^\star$ is given by
\begin{equation}
{\mathbb F}{\mathcal L}\left(x^0,x^1,\phi,
\frac{\partial \phi}{\partial x^0}, \frac{\partial \phi}{\partial x^1}\right)
 = \left[x^0,x^1,
p^0 \frac{\partial \phi}{\partial x^0}+p^1 \frac{\partial \phi}{\partial x^1}
L\left(x^0,x^1,\phi,\frac{\partial \phi}{\partial x^0}, \frac{\partial
\phi}{\partial x^1}\right) ,
\phi,\frac{\partial \phi}{\partial x^0}, \frac{\partial \phi}{\partial
x^1}\right]
\nonumber
\end{equation}
where $p^\mu= \partial L / \partial \phi_\mu$ and $\phi_\mu = \partial \phi
/ \partial
x^\mu$.  When $L$ is regular (it is for the nonlinear wave equation), we
have the
primary constraint subbundle
$P_{\mathcal L}:= \pi_{{\mathbb R}^2,{\mathbb R}^2 \times {\mathbb R}^3 }
 \subset J^1(Y)^\star$ with coordinates
$(x^0,x^1,\phi,p^0,p^1)$, and the Hamiltonian density on $P_{\mathcal L}$ is
written in coordinates as
\begin{equation}
{\mathcal H}(x^0,x^1,\phi,p,p^0,p^1) =
\left[p^0 \frac{\partial \phi}{\partial x^0}+p^1 \frac{\partial \phi}{\partial
x^1} - L\left(x^0,x^1,\phi,\frac{\partial \phi}{\partial x^0}, \frac{\partial
\phi}{\partial x^1}\right)\right]
 dx^1 \wedge x^0,
\nonumber
\end{equation}
while the canonical $3$-form on $J^1(Y)^\star$ is given by
\begin{equation}
\Omega_{\mathcal H}=-d\phi\wedge dp^0\wedge dx^1 + d\phi\wedge dp^1 \wedge dx^0
+ dH\wedge dx^1 \wedge dx^0 .\label{E5}
\end{equation}
We note that in this case, by global triviality, we may identify
$P_{\mathcal L}$ with
$\pi_{{\mathbb R}^2,{\mathbb R}^3 }$.

Bridges [1997] considers this scalar field theory with the manifold
$M={\mathbb R}$ and a Lagrangian $L$ that has no explicit dependence on time or
space.  (In particular, all of the fibers of both $J^1(Y)$ and $P_{\mathcal L}$
are identical over $X$ and identified with
${\mathbb R}^3$.)  He obtains the following partial differential
equation for ${\mathcal Z}\equiv (\phi, p^0,p^1)$ :
\begin{equation}
M\frac{\partial{\mathcal Z}}{\partial x^0} + K \frac{\partial {\mathcal
Z}}{\partial x^1}
=-dH ({\mathcal Z}), \label{E2}
\end{equation}
where the $3\times 3$ matrices $M$ and $K$ are defined to be
\begin{equation}
M =
\left[\begin{array}{cccccc}
0   & -1  & 0\\
1   & 0  & 0\\
0   & 0  & 0\\
\end{array} \right] , \ \ \
K =
\left[\begin{array}{cccccc}
0   & 0  & -1\\
0   & 0  & 0\\
1  & 0  & 0\\
\end{array} \right].
\label{E3}
\end{equation}
$M$ and $K$ may be identified with a pair of degenerate $2$-forms
$\omega^{(1)}$ and $\omega^{(2)}$ on
$P_{\mathcal L}$ which define Bridges' multisymplectic structure, and although
it may appear that these two $2$-forms provide a distinct structure
from that of the $3$-form in (\ref{E5}), in fact it is just a
particular coordinate representation of the intrinsic structure which we
have defined.
\begin{prop}
\label{PROP2}
If ${\mathbb F}{\mathcal L}:J^1(Y)\rightarrow P_{\mathcal L}$ is a fiber bundle
diffeomorphism and  $\phi \in \Gamma(\pi_{XY})$, then
$\tilde{j}^1(\phi)$ is a Hamiltonian system for ${\mathcal H}$ if and only if
\begin{equation}
\frac{\partial}{\partial x^0} \hook \tilde{j}^1(\phi)^* (dp^0\wedge d\phi)
+
\frac{\partial}{\partial x^1} \hook \tilde{j}^1(\phi)^* (dp^1\wedge d\phi)
= -dH(\tilde{j^1}(\phi)).
\label{E4}
\end{equation}
where (\ref{E4}) is equivalent to Bridges' equation (\ref{E2}).
\end{prop}
\begin{proof}
Let $U \in T(P_{\mathcal L})$ be an arbitrary vector field which in
coordinates is
\begin{equation}
U=U_{x^0}\frac{\partial}{\partial x^0} + U_{x^1}\frac{\partial}{\partial x^1}+
U_\phi
\frac{\partial}{\partial \phi} +U_{p^0}\frac{\partial}{\partial p^0}+
U_{p^1}\frac{\partial}{\partial p^1},
\nonumber
\end{equation}
so that
\begin{eqnarray*}
U\hook \Omega_{\mathcal H}&=&-U_\phi dp^0\wedge dx^1 + U_{p^0} d\phi\wedge dx^1
-U_{x^1} d\phi \wedge dp^0 \\
&&\  +\ U_\phi dp^1\wedge dx^0 - U_{p^1} d\phi\wedge dx^0
+U_{x^0} d\phi \wedge dp^1 \\
&&\  +\ \left(\frac{\partial H}{\partial \phi}  U_\phi
+\frac{\partial H}{\partial p^0}  U_{p^0}
+\frac{\partial H}{\partial p^1}  U_{p^1} \right) dx^1\wedge dx^0\\
&&\  -\ \frac{\partial H}{\partial \phi}  U_{x^1}d\phi\wedge dx^0
-\frac{\partial H}{\partial p^0}  U_{x^1} dp^0\wedge dx^0
-\frac{\partial H}{\partial p^1}  U_{x^1}  dp^1\wedge dx^0\\
&&\  +\ \frac{\partial H}{\partial \phi}  U_{x^0}d\phi\wedge dx^1
+\frac{\partial H}{\partial p^0}  U_{x^0} dp^0\wedge dx^1
+\frac{\partial H}{\partial p^1}  U_{x^0}  dp^1\wedge dx^1.
\end{eqnarray*}
Then
\begin{eqnarray*}
\tilde{j}^1(\phi)^*(U\hook \Omega_{\mathcal H}) &=&
\left[
U_\phi\left(\frac{\partial p^0}{\partial x^0}+\frac{\partial p^1}
{\partial x^1} +\frac{\partial H}{\partial \phi}\right) +
U_{p^0}\left(-\frac{\partial \phi}{\partial x^0}+ \frac{\partial H}
{\partial p^0}\right) + U_{p^1}\left(-\frac{\partial \phi}{\partial x^1}
+ \frac{\partial H} {\partial p^1}\right) \right.\\
&&+\ U_{x^0}\left( -\frac{\partial \phi}{\partial x^0}
\frac{\partial p^1}{\partial x^1}
+ \frac{\partial \phi}{\partial x^1} \frac{\partial p^1}{\partial x^0}
-  \frac{\partial H}{\partial \phi} \frac{\partial \phi}{\partial x^0}
-  \frac{\partial H}{\partial p^0} \frac{\partial p^0}{\partial x^0}
- \frac{\partial H}{\partial p^1} \frac{\partial p^1}{\partial x^0}
\right)\\
&&\ \left. +\ U_{x^1}\left( \frac{\partial \phi}{\partial x^0}
\frac{\partial p^0}{\partial x^1}
-\frac{\partial \phi}{\partial x^1} \frac{\partial p^0}{\partial x^0}
- \frac{\partial H}{\partial \phi} \frac{\partial \phi}{\partial x^1}
- \frac{\partial H}{\partial p^0} \frac{\partial p^0}{\partial x^1}
- \frac{\partial H}{\partial p^1} \frac{\partial p^1}{\partial x^1}
\right) \right]dx^1 \wedge dx^0,
\end{eqnarray*}
which vanishes if and only if
\begin{equation}
dH = -\left(\frac{\partial p^0}{\partial x^0}+\frac{\partial p^1}
{\partial x^1}\right) d\phi + \frac{\partial \phi}{\partial x^0} dp^0+
\frac{\partial \phi}{\partial x^1} dp^1, \nonumber
\end{equation}
and this is precisely a restatement of (\ref{E4}). To see that (\ref{E4}) is
equivalent to (\ref{E2}), simply notice that $\partial_\mu \tilde{j^1}(\phi)
= T \tilde{j^1}(\phi) \cdot \partial_\mu$ and that $\partial_\mu {\mathcal Z}$
is the $\pi_{X,P_{\mathcal L}}$-vertical component of
$\partial_\mu \tilde{j^1}(\phi)$.
\end{proof}
Thus, we have shown that Bridges' formulation is equivalent to our
intrinsically defined multihamiltonian system for the nonlinear
wave equation defined over one spatial dimension ($n=1$).  The
argument, however, is entirely independent of  the number of spatial
directions, and obviously holds when $X={\mathbb R}^{n+1}$ and
$Y={\mathbb R}^{n+1} \times {\mathbb R}$, in which
case our multihamiltonian system may be expressed as
$$\frac{\partial }{\partial x^\mu} \hook \tilde{j^1}(\phi)^*
(dp^\mu \wedge d\phi) = -dH,$$
or, in terms of Bridges' $n$$+$$1$ $2$-forms $\omega^{(\mu)}$, as
$$ \omega^{(\mu)} \left(\frac{\partial\tilde{j^1}(\phi)}{\partial x^\mu}
, U \right) = -dH(\tilde{j^1}(\phi)) \cdot U \text{ for all } U \in
T(P_{\mathcal L}).$$
More importantly, as we shall show in
Section \ref{sec6}, in the presence of symmetry, we can assemble
these $n$$+$$1$ distinct $2$-forms $\omega^{(\mu)}$ into
our single ($n$$+$$2$)-form $\Omega_{\mathcal H}$.

\section{Covariant Noether Theory}
\label{sec5}
\begin{defn}
A {\bf covariant canonical transformation} is a $\pi_{XZ}$-bundle map
$\eta_Z : Z \rightarrow Z$ covering a diffeomorphism $\eta_X : X \rightarrow X$
such  that $\eta^*_Z\Omega  =  \Omega$. $\Diamond$
\end{defn}
\begin{defn}
If $\eta_Y :  Y \rightarrow Y$ is a $\pi_{XY}$-bundle automorphism (also
covering a diffeomorphism $\eta_X : X \rightarrow X$), its {\bf canonical lift}
$\eta_Z : Z \rightarrow Z$ is defined by
\begin{equation}
\eta_Z(z)  =  (\eta_Y^{-1})^*(z). \qquad \Diamond \label{F1}
\end{equation}
\end{defn}

We may now define the covariant analogue of momentum maps in symplectic
geometry.

Let ${\mathcal G}$ denote a Lie group (perhaps infinite-dimensional) with Lie
algebra ${\mathfrak g}$  that acts on  $X$ by diffeomorphisms and acts on
$Z$ (or
$Y$) as $\pi_{XZ}$ (or $\pi_{XY})$-bundle automorphisms. For $\eta \in
{\mathcal G}$, let $\eta_X, \eta_Y$ and $\eta_Z$ denote the corresponding
transformations of $X, Y$ and $Z$ (the map $\eta_Z:Z \rightarrow Z$ is the
prolongation of $\eta_Y$) and for $\xi \in {\mathfrak g}$,  let $\xi_X,
\xi_Y$ and $\xi_Z$ denote the corresponding infinitesimal generators.
If ${\mathcal G}$ acts on $Z$ by covariant canonical transformations, then the
Lie derivative of $\Omega$  along $\xi_Z$ is zero:
\begin{equation}
{\mathfrak L}_{\xi_Z} \Omega = 0, \label{F2a}
\end{equation}
so that the left Lie algebra action is canonical.  In the case that
\begin{equation}
{\mathfrak L}_{\xi_Z} \Theta = 0, \label{F2b}
\end{equation}
then ${\mathcal G}$ acts by {\it special} covariant transformations.

\begin{defn}
Let a Lie algebra ${\mathfrak g}$ have a canonical left action on $Z$, and
suppose
there exists $J\in L({{\mathfrak g}},\Lambda^n (Z))$ covering the identity on
$Z$ such
that for each $\xi \in {\mathfrak g}$,
\begin{equation}
\xi_Z \hook \Omega = dJ(\xi ). \nonumber
\end{equation}
The map ${\mathbf J}:Z \rightarrow {\mathfrak g}^* \otimes \Lambda^n(Z)$
defined by
\begin{equation}
{\bold J}(z)(\xi )  = J(\xi ) (z)
\label{F3}
\end{equation}
for all $\xi \in {\mathfrak g}$ and $z\in Z$ is called a {\bf covariant
momentum mapping} (or a {\bf multimomentum mapping})  of the action.
$\Diamond$
\end{defn}
The covariant momentum map is said to be Ad$^*$-equivariant if the
diagram
\begin{equation}
\begin{array}{ccc}
Z & \stackrel{{\bf J}}{\longrightarrow} & {\mathfrak g}^*
\otimes\Lambda^{n}(Z) \\
\Big\downarrow\vcenter{%
  \rlap{$\scriptstyle{\rm }\,\eta_Z$}} &
        & \Big\downarrow\vcenter{%
           \rlap{$\scriptstyle{\rm }\,\text{Ad}^* \otimes \text{Id}$}}\\
Z & \stackrel{{\bf J}}{\longrightarrow} &{\mathfrak g}^* \otimes \Lambda^{n}(Z)
\end{array}
\nonumber
\end{equation}
commutes, or equivalently if
\begin{equation}
J(\text{Ad}_\eta^{-1} \xi) = \eta_Z^*[J(\xi)].  \label{F3b}
\end{equation}

\begin{lemma}
If the action on $Z$ is the lifted action $\eta_Z$, then ${\mathcal G}$ acts by
special covariant transformations, the mapping $J$ defined by
\begin{eqnarray}
J( \xi ) & = &\xi_Z \hook \Theta \nonumber
\\ & = & \pi_{YZ}^* ({\xi_Y} \hook z), \label{F12}
\end{eqnarray}
is a multimomentum mapping of the action for the multisymplectic form
on $\Omega$, i.e.,
\begin{equation}
\xi_Z \hook \Omega = d J(\xi),
\label{S0}
\end{equation}
and is Ad$^*$-equivariant.
\end{lemma}
\begin{proof}
Differentiating the coordinate expression for (\ref{F1}) we find that
if $\xi = (\xi^\mu,\xi^A)$, then
\begin{equation}
\xi_Z = (\xi^\mu,\xi^A, - p\xi^\nu{}_{,\nu} - p_B{}^\nu \xi^B{}_{,\nu},
p_A{}^\nu \xi^\mu{}_{,\nu} - p_B{}^\mu \xi^B{}_{,A} -
p_A{}^\mu \xi^\nu{}_{,\nu}), \label{F10}
\end{equation}
and hence that ${\mathfrak L}_{\xi_Z}\Theta=0$.  Then,
\begin{equation}
dJ(\xi ) = d({\xi_Z}\hook\Theta) = {\mathfrak
L}_{\xi_Z} \Theta -{\xi_Z} \hook d\Theta = {\xi_Z}\hook \Omega. \nonumber
\end{equation}
Since $\xi_Y = T\pi_{YZ}\circ \xi_Z$, (\ref{F12}) immediately follows,
and the last assertion holds because special covariant momentum maps
are $Ad^*$-equivariant (the argument is analogous to that for the
cotangent bundle case which is proven in Abraham and Marsden [1978],
Theorem 4.2.10).
\end{proof}

In coordinates this {\bf special covariant momentum map} may be
expressed as
\begin{equation}
J(\xi)(z) = ({p_A}^\mu \xi^A + p \xi^\mu)d^nx_\mu
- {p_A}^\mu
\xi^\nu dy^A \wedge {\partial_{\nu}}\hook ({\partial_{\mu}}\hook d^{n+1}x).
\label{F14}
\end{equation}
Next we describe the prolonged action of the group ${\mathcal G}$ on
$Y$ to $J^1(Y)$ and $P_{\mathcal L}$.
\begin{defn}
\label{prolong}
Let $\eta_Y:Y\rightarrow Y$ be a $\pi_{XY}$-bundle automorphism
covering a diffeomorphism $\eta_X:X\rightarrow X$.
Then
\begin{equation}
\eta_{J^1(Y)}(\gamma)=T\eta_Y \circ \gamma \circ T\eta_X^{-1}
\text{ for all } \gamma \in J^1(Y), \label{S1}
\end{equation}
and
\begin{equation}
\eta_{P_{\mathcal L}}(z)=
i_{J^1(Y)^\star,P_{\mathcal L}}^* {\mathbb F}{\mathcal L}
(\eta_{J^1(Y)}\circ
({\mathbb F}{\mathcal L}|_{P_{\mathcal L}})^{-1} z)
\text{ for all } z \in P_{\mathcal L}.  \Diamond\label{S2}
\end{equation}
\end{defn}

\begin{defn}
We say that the Lagrangian density ${\mathcal L}$ and the Hamiltonian
density ${\mathcal H}$ are
equivariant with respect to ${\mathcal G}$ if for all $\eta \in
{\mathcal G}$, $\gamma\in J^1(Y)$, and $z\in P_{\mathcal L}$,
\begin{equation}
{\mathcal L}(\eta_{J^1(Y)}(\gamma)) = (\eta_{X}^{-1})^* {\mathcal L}(\gamma),
\label{S3}
\end{equation}
and
\begin{equation}
{\mathcal H}(\eta_{P_{\mathcal L}}(z)) = (\eta_{X}^{-1})^* {\mathcal H}(z),
\label{S4}
\end{equation}
where $(\eta_{X}^{-1})^*{\mathcal L}(\gamma)$ means the $(n+1)$-form ${\mathcal
L}(\gamma)$
at $x \in X$ is pushed forward to an $(n+1)$-form at $\eta(x)$. $\Diamond$
\end{defn}
Analogous to Corollary 4.2.14 of Abraham and Marsden [1978], one may readily
verify that both $\Theta_{\mathcal L}$ and $\Theta_{\mathcal H}$ are invariant
under the respective  group action prolongations, i.e.,
\begin{equation}
\eta_{J^1(Y)}^* \Theta_{\mathcal L} = \Theta_{\mathcal L} \text{ and }
\eta_{P_{\mathcal L}}^* \Theta_{\mathcal H}=\Theta_{\mathcal H}. \nonumber
\end{equation}

\begin{lemma}
\label{equi}
Suppose that ${\mathcal L}$ is regular and that ${\mathcal H}:P_{\mathcal L}
\rightarrow \Lambda^{n+1}(X)$ is equivariant, and is not constant on any
$P_{\mathcal L}$ neighborhood.  Then ${\mathcal L}:J^1(Y)\rightarrow
\Lambda^{n+1}(X)$ is equivariant.
\end{lemma}
\begin{proof}
Let ${\mathcal G}$ be a group acting on $Y$ by bundle automorphisms.  From
(\ref{S2}),
we have by definition that  ${\mathbb F}{\mathcal L}:J^1(Y)\rightarrow
P_{\mathcal L}$
is equivariant with respect to ${\mathcal G}$.  Hence, for all $\eta \in
{\mathcal G}$
and $\gamma \in J^(Y)$,
\begin{equation}
{\mathcal H}(\eta_{P_{\mathcal L}}({\mathbb F}{\mathcal L}(\gamma))) =
(\eta_{X}^{-1})^* {\mathcal H}({\mathbb F}{\mathcal L}(\gamma)).
\label{S5}
\end{equation}
Assume that the Lagrangian density ${\mathcal L}$ is not equivariant with
respect to
${\mathcal G}$.  Then, there exists $\eta \in {\mathcal G}$ and $\gamma \in
J^1(Y)$
for which
\begin{equation}
F (\gamma) \equiv \left[
{\mathcal L} \circ \eta_{J^1(Y)}(\gamma) -(\eta_{X}^{-1})^* {\mathcal L}
\right] (\gamma)
\neq 0.  \nonumber
\end{equation}
Hence, by continuity of $F$, there is some neighborhood $U$ in $J^1(Y)$
about $\gamma$
for which $F(U)$ does not intersect $\{0\}$ in $\Lambda^{n+1}(X)$.  We will
assume that
$\gamma \in J^1(Y)_y$ for some fixed $y\in Y$ and take $U$ to be a fiber
neighborhood
of $\gamma$.  By choosing $\varepsilon > 0$ sufficiently small, we see that
for all
$\gamma' \in J^1(Y)_y$ satisfying ${\mathfrak A} \circ \gamma' = {\mathfrak
A}$,
\begin{equation}
\left. \frac{d}{d\varepsilon}\right|_{\varepsilon =0}
{\mathcal L} (\eta_{J^1(Y)}(\gamma) + \varepsilon(\eta_{J^1(Y)}(\gamma')
-\eta_{J^1(Y)}(\gamma)) \neq
(\eta_{X}^{-1})^* \left. \frac{d}{d\varepsilon}\right|_{\varepsilon =0}
{\mathcal L} (\gamma + \varepsilon(\gamma' -\gamma)), \nonumber
\end{equation}
in which case,
$$ D_{\mathfrak A} (\eta_{J^1(Y)}(\gamma)) \cdot \eta_{J^1(Y)}(\gamma) \neq
(\eta_X^{-1})^* D_{\mathfrak A} (\gamma)\cdot \gamma $$
for all $\gamma$ in some $J^1(Y)$-neighborhood $U$.  Since ${\mathcal H}$
cannot be locally constant, (\ref{S5}) cannot be true for all $\gamma \in
U$, and
so ${\mathcal H}$ cannot be equivariant.
\end{proof}

\begin{thm}
\label{thm_mom}
Let ${\mathcal L}$ be a regular Lagrangian density, and suppose that
${\mathcal H}:P_{\mathcal L}\rightarrow\Lambda^{n+1}(X)$ is equivariant
and is not constant on any $P_{\mathcal L}$-neighborhood.  Then the
map
$$J^{\mathcal L}(\xi):={\mathbb F}{\mathcal L}^*J(\xi)$$
is a momentum map for the lifted action of ${\mathcal G}$ on $J^1(Y)$
relative to $\Omega_{\mathcal L}$, and  the map
$$J^{\mathcal H}( \xi ) :=(({\mathbb F}{\mathcal L}
|_{P_{\mathcal L}})^{-1})^* \circ {\mathbb F}{\mathcal L}^* J(\xi)
=({\mathbb F}{\mathcal L}\circ {\mathbb F}{\mathcal L}|_{P_{\mathcal L}}
^{-1})^* J(\xi)= i_{J^1(Y)^\star,P_{\mathcal L}}^* J(\xi)$$
is a momentum map for the lifted action of ${\mathcal G}$ on $P_{\mathcal L}$
relative to $\Omega_{\mathcal H}$; i.e., for all $\xi \in {\mathfrak g}$,
\begin{equation}
\xi_{J^1(Y)} \hook \Omega_{\mathcal L} = d J^{\mathcal L}(\xi),
\label{S6}
\end{equation}
and
\begin{equation}
\xi_{P_{\mathcal L}} \hook \Omega_{\mathcal H} = d J^{\mathcal H}(\xi),
\label{S7}
\end{equation}
where $\xi_{J^1(Y)}$ and
where $\xi_{P_{\mathcal L}}$ are the infinitesimal generators
corresponding to $\xi$.
\end{thm}
\begin{proof}
Lemma \ref{equi} asserts that ${\mathcal L}$ is equivariant, from
which we may conclude that  ${\mathbb F}{\mathcal L}:J^1(Y)
\rightarrow Z$ is equivariant so that
\begin{equation}
\eta_Z \circ {\mathbb F}{\mathcal L} = {\mathbb F}{\mathcal L} \circ
\eta_{J^1(Y)}. \label{N7}
\end{equation}
Indeed, we see that
\begin{eqnarray*}
\{\eta_{J^1(Y)^\star}[{\mathbb F}{\mathcal
L}(\gamma)]\}\cdot\gamma'&=&(\eta^{-1}_X)^*
\{{\mathbb F}{\mathcal L}(\gamma)[\eta^{-1}_{J^1(Y)}(\gamma')]\} \\
&=& (\eta^{-1}_X)^*\left\{{\mathcal L}(\gamma) + \left.\frac{d}{d\varepsilon}
\right|_{\varepsilon =0}
{\mathcal L}(\gamma + \varepsilon[\eta^{-1}_{J^1(Y)}(\gamma') - \gamma])
\right\}.
\end{eqnarray*}
and that
\begin{equation}
\{{\mathbb F}{\mathcal L}[\eta_{J^1(Y)}(\gamma)]\}\cdot\gamma' =
 {\mathcal L} (\eta_{J^1(Y)}(\gamma)) + \left. \frac{d}{d\varepsilon}\right|
_{\varepsilon=0}
 {\mathcal L}(\eta_{J^1(Y)}(\gamma)) + \varepsilon[\gamma' - \eta_{J^1(Y)}
(\gamma)]) ,
\end{equation}
which are equal by the equivariance of ${\mathcal L}$.  The
infinitesimal version of (\ref{N7}) yields
\begin{equation}
\xi_Z \circ {\mathbb F}{\mathcal L} = T{\mathbb F}{\mathcal L} \cdot
\xi_{J^1(Y)},
\label{N8}
\end{equation}
which is a statement that $\xi_Z$ is ${\mathbb F}{\mathcal L}$-related
to $\xi_{J^1(Y)}$.  Hence, the pull-back of (\ref{S0}) along
${\mathbb F}{\mathcal L}$ gives us (\ref{S6}), while a second
pull-back of (\ref{S6}) along the diffeomorphism
${\mathbb F}{\mathcal L}|_{P_{\mathcal L}}^{-1}$ verifies (\ref{S7}).
\end{proof}

In multisymplectic coordinates, the multimomentum mapping $J^{\mathcal H}$
is written as
\begin{equation}
J^{\mathcal H}(\xi) = ({p_A}^\mu \xi^A + ({p_A}^\nu {\mathfrak A}^A_\nu
-H) \xi^\mu)d^nx_\mu - {p_A}^\mu
\xi^\nu dy^A \wedge {\partial_{\nu}}\hook ({\partial_{\mu}}\hook d^{n+1}x),
\label{N10}
\end{equation}
where $H:P_{\mathcal L}\rightarrow {\mathbb R}$ is the Hamiltonian function
associated with the Hamiltonian density ${\mathcal H}$ by ${\mathcal H}
=H d^{n+1}x$.

\begin{thm}[Covariant Hamiltonian Version of Noether's Theorem]
\label{noether}
Let ${\mathcal L}$ be regular, and suppose that
a group ${\mathcal G}$ acts on $Y$ by $\pi_{XY}$-bundle automorphisms,
and that the
Hamiltonian density ${\mathcal H}$ is equivariant with respect to
${\mathcal G}$, and is not locally constant on any
$P_{\mathcal L}$-neighborhood.  Then for each $\xi \in {\mathcal G}$,
\begin{equation}
d[\tilde{j}^1(\phi)^* J^{\mathcal H} (\xi )]=0,
\label{F15}
\end{equation}
for any $\phi \in \Gamma(\pi_{XY})$ for which $\tilde{j}^1(\phi)$ is a
covariant  Hamiltonian system.
\end{thm}
\begin{proof}
Since $\tilde{j}^1(\phi)$ is a Hamiltonian system, $\tilde{j}^1(\phi)^* i_U
\Omega_{\mathcal H} =0$ for any vector $U \in T(P_{\mathcal L})$ so set
$U=\xi_{P_{\mathcal L}}.$
\end{proof}
Under the same hypotheses, Lemma \ref{LEMMA1} gives us the following equivalent
statement.
\begin{cor}
\label{cor5.1}
For each $\xi \in {\mathcal G}$, $d[j^1(\phi)^* J^{\mathcal L}(\xi)]=0$.
\end{cor}
The quantity $\tilde{j^1}(\phi)^* J^{\mathcal H}(\xi)$ is called the
{\bf Hamiltonian Noether current}, and as we shall show,
leads to very useful decompositions of the classical water wave
conservations laws.  In coordinates, it has the form
\begin{equation}
\tilde{j^1}(\phi)^* J^{\mathcal H}(\xi)=
\left[
({p_A}^\mu\xi^A +({p_A}^\nu {\mathfrak A}^A_\nu -H)\xi^\mu) -
{p_A}^\mu {\phi^A}_{,\nu} \xi^\nu +
{p_A}^\nu {\phi^A}_{,\nu} \xi^\mu \right]d^nx_\mu. \label{S10}
\end{equation}
\section{Symmetry and generalized conservation laws}
\label{sec6}
Just as we have shown in Section \ref{sec4} that the
multihamiltonian system generalizes the classical
Hamiltonian description of particle mechanics as well as the
structures defined in Bridges [1997], we can do the same for the
multisymplectic Hamiltonian Noether theory.
\subsection{Particle mechanics}
\label{subsec6.1}
We let the groups Diff$({\mathbb R})$ act on ${\mathbb R}$ and $G$ on
$Q$ and consider the action of the prolongation of ${\mathcal G} =$
Diff$({\mathbb R})
\times G$.  With elements of ${\mathfrak g}$ written as
$(f,\xi)$,
the Hamiltonian Noether current has the simple coordinate form
\begin{equation}
\tilde{j^1}(\phi)^* J^{\mathcal H}(\xi)= p_A \xi^A - Hf := J^H_\xi  - Hf,
\nonumber
\end{equation}
where $J^H_\xi = p_A \xi^A$ is the usual momentum map for $G$
acting on $Q$ in Hamiltonian mechanics.  Then equation (\ref{F15}) asserts
that along  trajectories of the Hamiltonian vector field,
\begin{equation}
\frac{d}{dt}( J^H_\xi -Hf)=0, \nonumber
\end{equation}
which is equivalent to conservation of {\it both} $J^H_\xi$ and $H$.

\subsection{Bridges' decomposition of Noether theory}
\label{subsec6.2}
In this section we fix
$X={\mathbb R}^{n+1}$, $Y={\mathbb R}^{n+1} \times {\mathbb R}$, and
$P_{\mathcal L} = \pi_{{\mathbb R}^{n+1},{\mathbb R}^{n+2}}$.  We examined
this particular geometry in Section \ref{subsec_pde}, wherein we showed that Bridges' 
multisymplectic structure consisting of $n$$+$$1$ degenerate $2$-forms $\omega^{(\mu)}$
on $P_{\mathcal L}$ is an equivalent and particular 
coordinate representation of our intrinsically defined multihamiltonian
system on $P_{\mathcal L}$.  We  now show that 
in the presence of a general symmetry group action, 
Bridges' $n$$+$$1$ pre-symplectic $2$-forms $\omega^{(\mu)}$ can actually be assembled
into our single multisymplectic $(n+2)$-form $\Omega^{\mathcal H}$.

Bridges forms an $(n+1)$-vector of momentum mappings such that each vector 
component is associated with a distinct spacetime direction through the 
pre-symplectic $2$-form $\omega^{(\mu)}$.  He then relates the action of
a general symmetry group ${\mathcal G}$ along the fiber of $P_{\mathcal L}$
with the vanishing of the divergence of the vector of momentum mappings. The
following result is Theorem 2.2 in Bridges [1996b].
\begin{prop}
\label{prop_bridges}
Let $H:{\mathbb R}^{n+2} \rightarrow {\mathbb R}$ be a covariant Hamiltonian
with $n$$+$$1$ distinct $2$-forms $\omega^{(\mu)}$.  If $dH \cdot
\xi_{P_{\mathcal L}} = 0$ for all $\xi_{P_{\mathcal L}}$ in the Lie algebra
${\mathfrak g}$ of the group ${\mathcal G}$ acting on
${\mathbb R}^{n+2}$, and if $P^\mu$ is the momentum mapping associated to
$\omega^{(\mu)}$, i.e., for all $\xi_{P_{\mathcal L}} \in {\mathfrak g}$,
\begin{equation}
\xi_{P_{\mathcal L}} \hook \omega^{(\mu)} = dP^\mu(\xi_{P_{\mathcal L}}),
\label{Bmom}
\end{equation}
then
\begin{equation}
\frac{\partial P^\mu}{\partial x^\mu} = 0.
\label{S20}
\end{equation}
\end{prop}

We will first show that (\ref{S20}) is a particular example of our
conservation law (\ref{F15}) in the case of the trivial bundle geometry
defined above, for which the fields have no
explicit dependence on time or space (geometrically, this means that
each fiber of the bundle $\pi_{X,P_{\mathcal L}}$ is identical).
\begin{prop}
\label{prop5.2}
In the case that the action on $P_{\mathcal L}:=
\pi_{ {\mathbb R}^{n+1},{\mathbb R}^{n+2}}$ is
the lifted action $\eta_{P_{\mathcal L}}$, then the momentum
mappings $P^\mu$ defined in (\ref{Bmom}) are the components of the Hamiltonian
Noether current, and hence the conservation law (\ref{S20}) is contained
in Theorem \ref{noether}.
\end{prop}
\begin{proof}
We set the diffeomorphism $\eta_{X}$ to be the identity, and identify
$P_{\mathcal L}$ with ${\mathbb R}^{n+2}$.
Hence, the infinitesimal generators $\xi^\mu=0$, and
using (\ref{S10}), we see
that the Hamiltonian Noether current is given by
\begin{equation}
\tilde{j^1}(\phi)^* J^{\mathcal H}(\xi)= p^\mu \xi d^nx_\mu
:= N^\mu d^nx_\mu.
\label{Z2}
\end{equation}
Using (\ref{F10}), we easily deduce that
the lifted action $\xi_{P_{\mathcal L}}$ is given in coordinates by
$(0,\xi, -p^\mu \frac{\partial \xi}{\partial \phi})$ so that the equivariance
of ${\mathcal H}$ is equivalent to $dH \cdot \xi_{P_{\mathcal L}}=0$.

In accordance with
Proposition \ref{prop_bridges}, all of the group action is along the fiber of
$\pi_{ X, P_{\mathcal L}}$, identified with ${\mathbb R}^{n+2}$, 
so we will restrict the exterior derivative $d$ to the fiber.

We claim that $P^\mu = N^\mu$.  To see this, we must show that
$dN^\mu = \xi_{P_{\mathcal L}} \hook \omega^{(\mu)}$,
but this is precisely the case since in coordinates, for each $\mu=
1,...,n,0$,
$$ d N^\mu = \left( p^\mu \frac{\partial \xi}{\partial \phi}, 0, ..., \xi,
..., 0 \right),  \text{   $\xi$ in the $(\mu +1)$th coordinate}.$$ 
Then, using the identity
$$ d(N^\mu d^nx_\mu) = \partial_\nu N^\mu dx^\nu \wedge d^nx_\mu = \partial_\mu
N^\mu d^{n+1}x$$
we have that $d[\tilde{j^1}(\phi)^* J^{\mathcal H}(\xi)]=0$ implies
that $\frac{\partial N^\mu}{\partial x^\mu} = 0$ so that
$\frac{\partial P^\mu}{\partial x^\mu} = 0$ and the result is proved.
\end{proof}

This proposition indicates how we can assemble the $n$$+$$1$
$2$-forms $\omega^{(\mu)}$ into the single ($n$$+$$2$)-form
$\Omega_{\mathcal H}$ when lifted symmetries exist.  Namely,
to each $\omega^{(\mu)}$, there corresponds a momentum mapping
$P^\mu(\xi_{P_{\mathcal L}})$ of the symmetry group given
by (\ref{Bmom}).  By Proposition \ref{prop5.2}, the maps
$P^\mu(\xi_{P_{\mathcal L}})$ are the
components of the Hamiltonian Noether current.  This then
defines the Hamiltonian covariant momentum mapping $J^{\mathcal H}(\xi_
{P_{\mathcal L}})$ which in turn, by (\ref{S7}), defines the canonical
multisymplectic
($n$$+$$2$)-from $\Omega_{\mathcal H}$ on $P_{\mathcal L}$.  In fact,
since lifts are special canonical transformations, the covariant
momentum map defines the ($n$$+$$1$)-$\Theta_{\mathcal H}$ on
$P_{\mathcal L}$ as well.  We summarize with the following corollary.

\begin{cor}
\label{cor6.1}
Assume the group ${\mathcal G}$ acts by special canonical transformations, and
let the Hamiltonian Noether current $\tilde{j^1}(\phi)^* J^{\mathcal H}(\xi)$
be given in multisymplectic coordinates by $N^\mu d^n x_\mu$.
Then the $\omega^{(\mu)}$ satisfy
$\xi_{P_{\mathcal L}}\hook \omega^{(\mu)}= dN^\mu$.  Furthermore, if
$\omega^{(\mu)}$ is exact, such that
\begin{equation}
\omega^{(\mu)}= d \kappa^{(\mu)} \label{Z9}
\end{equation}
for  $1$-forms $\kappa^{(\mu)}$, then
\begin{equation}
N^\mu = \xi_{P_{\mathcal L}}\hook \kappa^{(\mu)}. \label{Z10}
\end{equation}
\end{cor}

\section{The Geometry of Water Waves}
\label{sec7}
It is interesting to note that the covariant Hamiltonian Noether
theorem intrinsically contains the mass conservation law for
water waves as well as the conservation of wave action and action
flux.  In particular, the vanishing of the exterior derivative
of the Hamiltonian Noether current is an intrinsic restatement of
the mass conservation
law, while the projected components of the Hamiltonian Noether
current $P^\mu$, as defined by Proposition \ref{prop5.2},  are
related to the action  and action flux.
As an example, for the case of two spatial dimensions,
the ensemble (or phase) average of $P^0$ corresponds to Whitham's
definition of wave action and that of $P^1$ and $P^2$ correspond
to the two-component action-flux (see Whitham [1974]), while in the case
of one spatial dimension the Hamiltonian density ${\mathcal H}$ is
related to the flow force or in some cases the momentum flux
(see Bridges [1997]).
These observations were first made by Bridges [1996a] in coordinates
and seemed to have been the primary motivating factors for defining
additional $2$-forms $\omega^{(\mu)}$ for each unbounded spatial
direction.

Next, we show that our definition for a multihamiltonian system
contains the variational principles which are essential to the
study of pattern formation and wave instability.  For simplicity,
we restrict our attention to symmetries given by the circle 
and ${\mathbb T}^{n+1}$; however,  it is important to note that our
procedure is general and may be applied to {\it any subgroup} of
the Euclidean group SE$(n+1)$ and its products.  This is significant
if one wishes to study hexagonal pattern formations, for example, in addition
to merely the periodic ones.

\subsection{Pattern formation, action, index, and the loop space}
\label{subsec7.1}
Let $X = {\mathbb R}^{n+1}$ and let $Y$ be the vector bundle ${\mathbb R}$ over
$X$.  Consider the semilinear elliptic scalar partial differential equation
\begin{equation}
\triangle \phi + V'(\phi) = 0, \ \ \phi \in \Gamma(\pi_{XY}),
\label{pde}
\end{equation}
where as above, $\triangle$ is the Laplace-Beltrami operator and
$V$ is a real-valued $C^\infty$-bundle map.  For this example,
it is appropriate to set
$P_{\mathcal L} = \pi_{{\mathbb R}^{n+1},{\mathbb R}^{n+2}}$ and
${\mathcal G}= SO(n+1)$, in which
case (\ref{pde}) may be equivalently expressed as
\begin{equation}
\tilde{j^1}(\phi)^* (U \hook \Omega_{\mathcal H}) = 0
\label{mspde}
\end{equation}
for all  $ U\in T(P_{\mathcal L})$, where in coordinates,
\begin{equation}
{\mathcal H} = Hd^{n+1}x \quad \text{  and  } \quad H(\tilde{j^1}(\phi))
= \frac{1}{2} p^\mu \cdot p^\mu + V(\phi).
\end{equation}
We show in this section that our intrinsic multisymplectic structure
can be used to generalize the notion of action and index on the
loop space of the primary constraint subbundle $P_{\mathcal L}$ as
defined in Bridges [1996b].

Let the map $\chi:X \rightarrow {\mathbb R}$ be defined in coordinates
by $\chi (x^\mu) = k_\mu x^\mu$, and identify ${\mathbb T}^1$ with its
universal cover ${\mathbb R}\backslash {\mathbb Z}$,  so that a smooth
$2\pi$-periodic map $\alpha: {\mathbb R}\rightarrow
P_{\mathcal L}$  may be identified with the smooth map $\alpha:
{\mathbb T}^1 \rightarrow  P_{\mathcal L}$.
\begin{defn}
The loop space of $P_{\mathcal L}$ is the subset of
$\Gamma(\pi_{X,P_{\mathcal L}})$
defined by
$$\text{Loop}(P_{\mathcal L}) = \{z \in \Gamma( \pi_{X,P_{\mathcal L} }) \mid
z= \alpha \circ \chi\text{ is holonomic and }
 \ \ \alpha \in C^\infty( {\mathbb T}^1, P_{\mathcal L} ) \}.$$
We then set
$$\text{loop}(P_{\mathcal L}) = \{z \in \text{Loop}(P_{\mathcal L}) \mid
\dot{\alpha}
\in {\mathcal P} \},$$
where $\dot{\alpha}= \frac{d\alpha}{d\chi}$ and
${\mathcal P} \subset T(P_{\mathcal L})$ consists of those vector fields on
$P_{\mathcal L}$ which are prolongations of vector fields on $Y$.
$\Diamond$
\end{defn}
\noindent
Hence, an element $\alpha \circ \chi$ in loop($P_{\mathcal L}$) is conjugate
to the first jet of a section $f \circ \chi$ in $\pi_{XY}$, where
$f:{\mathbb T}^1 \rightarrow Y$.

The diagonal periodic
patterns of (\ref{pde}) correspond to the restriction of (\ref{pde}) to
Loop$(P_{\mathcal L})$.   Thus, if $\alpha \circ \chi \in
\text{Loop}(P_{\mathcal L})$,
a diagonal periodic pattern satisfies
\begin{equation}
(\alpha \circ \chi)^* (U \hook \Omega_{\mathcal H}) = 0 \text{ for all }
U \in T(P_{\mathcal L}).
\label{loopPDE}
\end{equation}

Recently, the existence of periodic pattern solutions to (\ref{pde})
has been obtained by expressing such solutions as critical points of a
constrained variational principle, and using information provided by
sensitivity matrices, sometimes called the index, for classification
of the critical point type.  As it turns out, when the infinitesimal
group action coincides with the vector field $\dot{\alpha}$, the
Hamiltonian Noether current naturally and intrinsically verifies these
variational principles.

When $\alpha \circ \chi \in \text{loop}(P_{\mathcal L})$, we may
associate to it the loop space Hamiltonian Noether current
\begin{equation}
{\mathcal N}(\dot{\alpha}) := (\alpha \circ \chi)^* J^{\mathcal
H}(\dot{\alpha}).
\label{loopN}
\end{equation}
As with our previous example in Section \ref{subsec_pde}, we set the
group action on $X$ to be the identity, and identify $P_{\mathcal L}$
with ${\mathbb R}^{n+2}$.  We assume that ${\mathcal H}$ is equivariant
with respect to the group action on $P_{\mathcal L}$ and is not
locally constant.  In this case, $d {\mathcal N}(\dot{\alpha})=0$
implies that $\alpha \circ \chi := \tilde{j^1}(f \circ \chi)$ is a
Hamiltonian system, a fact that is readily verified by a tedious
computation in local coordinates. The computation, however, proceeds
easily on the Lagrangian side.   Using the coordinate expression
$$\frac{d}{d\chi}j^1(f \circ \chi)= j^1(\dot{f}):=\left(0, \dot{f},
\frac{\partial \dot{f}}{\partial x^\mu} +
\frac{\partial \dot{f}}{\partial y} \frac{\partial (f\circ \chi)}
{\partial x^\mu}\right), $$
which is readily obtained from the definition of the vector field
prolongation to $J^1(Y)$ given in (\ref{S1}), we appeal to
Corollary \ref{cor5.1} and check that
$d[j^1(f \circ \chi)^* J^{\mathcal L}(j^1(\dot{f}))]=0$.  In coordinates,
this yields
$$
\aligned
& \left\{ \left[\frac{\partial L}{\partial y} (j^1(f \circ \chi ))
- \frac{\partial}{\partial x^\mu}
\left( \frac{\partial L}{\partial (f\circ\chi)_\mu} (j^1(f \circ \chi))\right)
 \right] [-\dot{f} \circ (f \circ \chi)] \right. \\
& \qquad + \frac{\partial L}{\partial y}
(j^1(f \circ \chi)) (\dot{f} \circ (f \circ \chi)) \\
& \qquad + \left. \frac{\partial L}{\partial (f\circ\chi)_\mu} (j^1(f \circ
\chi))
\left[ \frac{\partial \dot{f}}{\partial x^\mu}
+ \frac{\partial \dot{f}}{\partial y} (f\circ\chi)_\mu
\right] (j^1(f \circ \chi))\right\} d^{n+1}x
\endaligned
$$
which vanishes if and only if the Lagrangian is equivariant and if
$f \circ \chi$ is a stationary point of $\int_X {\mathcal L}
(j^1(f\circ \chi) )$.  Since ${\mathcal H}$ is equivariant,
Lemma \ref{equi} guarantees that ${\mathcal L}$ is equivariant as well,
in which case Theorem \ref{THM1} gives us that $\alpha \circ \chi$ is
a Hamiltonian system.

To see that our covariant Noether theory contains the classical
constrained variational principle in local coordinates, we make the
following observations.  Let the $2$-forms $\omega^{(\mu)}$ and the
$1$-forms $\kappa^{(\mu)}$ be as defined in Corollary \ref{cor6.1}. By
Proposition \ref{prop_bridges}, equation (\ref{loopPDE}) is satisfied
if and only if
\begin{eqnarray*}
0&=& -dH(\alpha\circ\chi)  - k_\mu (\dot{\alpha}\hook \omega^{(\mu)}) \\
&=& d[-H(\alpha\circ\chi)  - k_\mu (\dot{\alpha}\hook \kappa^{(\mu)})]\\
&:=& d {\mathcal F}(\dot{\alpha}, \chi).
\end{eqnarray*}
Thus, as noted in Bridges [1996b], $\alpha\circ\chi$ is a diagonal periodic
pattern if it is the critical point of $\int_{{\mathbb T}^1}
{\mathcal F}(\dot{\alpha}, \chi) d\chi$ (classically, the phase-averaged
quantities are considered).  From this, we see that the solutions to
(\ref{loopPDE}) are the critical points of the phase-averaged Hamiltonian
with the additional constraints
\begin{equation}
\int_{{\mathbb T}^1} [\dot{\alpha}\hook \kappa^{(\mu)}] d \chi = I_\mu
\label{Z20}
\end{equation}
so that the $k_\mu$ are the Lagrange multipliers.  In the case that
$H$ may be viewed as an implicit function of the $I_\mu$, we have that
$$ k_\mu = -\frac{\partial H}{\partial I_\mu},$$
so that the Hessian matrix of H with respect to the level sets $I_\mu$
has components
$$ [\text{Hess}_I(H)]_{\mu \nu}= \frac{\partial^2 H}{\partial I_\mu
 \partial I_\nu} = \frac{\partial k_\mu }{\partial I_\nu}.$$
From the implicit function theorem,  a diagonal periodic pattern
is non-degenerate if det$[$Hess$_I(H)] \neq 0$.  Then, the natural
definition for the index for such patterns is given by
$$\text{index}(\alpha) = \#\text{ negative eigenvalues of Hess}_I(H).$$
We refer the reader to Bridges [1996a] for a detailed account and
applications.

\subsection{Stability of Water Waves}
\label{subsec7.2}
Conservative partial differential equations are often accurate models
for water waves, and in this section we will briefly comment upon
the connection between our covariant Noether theory and the
constrained toral variational principles which lead to characterizations
of the instabilities of the system.  Our brevity is due to the
fact that the  Hess$_I(H)$-matrix is explicitly connected to the
linear stability exponents from which we may deduce the behavior
of our solutions, and as we gave a fairly detailed description in the
previous section of how this matrix arises from the vanishing of the
exterior derivative of the Hamiltonian Noether current, we shall herein
only discuss the minor modifications necessary for this theory.

To demonstrate the main ideas, let us consider the
the manifolds $X$, $Y$, and $P_{\mathcal L}$ as given in the previous
section and the the partial differential equation defined in (\ref{Z1}),
with corresponding covariant Hamiltonian
$${\mathcal H} = \left[ \frac{1}{2}\left(\sum_{\mu=1}^n p_\mu^2  - p_0^2
 \right) - V(\phi) \right] d^{n+1}x.$$
Unlike the case of pattern formation for which we considered solutions
of
\begin{equation}
\tilde{j^1}(\phi)^*(U\hook \Omega_{\mathcal H})=0 \text{ for all } U
\in T(P_{\mathcal L}) \label{Z40}
\end{equation}
restricted to loop$(P_{\mathcal L})$, now we restrict consideration to
the periodic
sections of $P_{\mathcal L}$, so that $\tilde{j^1}(\phi):{\mathbb T}^{n+1}
\rightarrow P_{\mathcal L}$.  If we make the change of variables
$w^\mu=k_\mu x^\mu$ (no sum), then by Proposition \ref{prop_bridges},
periodic solutions of (\ref{Z40}) are expressed in coordinates by
\begin{equation}
k_\mu\frac{ \partial \tilde{j^1}(\phi)}{\partial w^\mu} \hook
\omega^{(\mu)} = -dH(\tilde{j^1}(\phi)), \ \ (k_0,...,k_n)\in
{\mathbb T}^{n+1}.
\label{Z50}
\end{equation}
Arguing exactly as we did in the previous section, we may again deduce
that the $k_\mu$ are the Lagrange multipliers of the system, and thus
the Hessian matrix of $H$ with respect to the level set $I_\mu$  is
identically obtained as for the case of periodic pattern formation.
See Bridges [1997] for a discussion of the relationship between
Hess$_I(H)$ and the classical linear stability exponents.
\vspace{.2in}

\noindent{\bf Acknowledgements.}
The authors would like to thank Tom Bridges and Mark Gotay for
reading early drafts of the manuscript and making numerous
comments and suggestions.
S.S. was partially supported by a fellowship
from the Cecil H. and Ida M. Green Foundation.
J.E.M. was partially supported by the National Science
Foundation under Grant DMS-DMS-96 33161 and the Department of
Energy under Contract DE-FG0395-ER25251.


\begin{thebibliography}{steve}
\bibitem{AM} {\sc Abraham, R. and J. Marsden}, [1978], {\em
Foundations of Mechanics}, Second Edition, Addison-Wesley, Menlo Park,
CA.
\bibitem{Bridges2} {\sc Bridges, T.J.}, [1996a], {\em
Periodic patterns, linear instability, symplectic structure
   and mean-flow dynamics for three-dimensional surface waves},
   Phil. Trans. R. Soc. Lond. A, 354, 533-574.
\bibitem{Bridges3} {\sc Bridges, T.J.}, [1996b], {\em
Symplecity, reversibility and elliptic operators}, Progress
   in Nonlinear Diff. Eq., 19, 1-20.
   \bibitem{Bridges1} {\sc Bridges, T.J.}, [1997], {\em
Multi-symplectic structures and wave propagation},
Math. Proc. Camb. Phil. Soc.\/, {\bf 121}, 147--190.
\bibitem{Bridges4} {\sc Bridges, T.J.}, {\em geometric formulation of
the conservation of wave action and its implications for signature
and the classification of instabilities},
to appear in Proc. R. Soc. Lond. A
\bibitem{CIL} {\sc Cantrijn, F., L.A. Ibort, and M. de Le\'{o}n}, 
[1997], {\em Hamiltonian  structures on multisymplectic manifolds},
to appear in Rend. Sem. Mat. Torino.
\bibitem{Gotay1} {\sc Gotay, M.J.}, [1991],
{\em A multisymplectic  framework for classical
field theory and the calculus of variations I: Covariant
Hamiltonian formulation}, Mechanics, Analysis and Geometry:
200 Years After Lagrange, M. Francaviglia, Ed., 203-235
(North Holland).
\bibitem{GIMMSY} {\sc Gotay, M., J. Isenberg, J.E. Marsden,
and R. Montgomery}, [1992], {\em Momentum Maps and the Hamiltonian Structure of
Classical Relativistic Field Theories.\/} (unpublished).
\bibitem{GM} {\sc Gotay, M.J. and J.E. Marsden}, [1992],
{\em Stress-energy-momentum tensors and the
Belifante-Resenfeld formula},
Cont. Math. AMS.\/ {\bf 132}, 367--392.
\bibitem{Whitham} {\sc Whitham, G.B.}, [1974] , {\em Linear and Nonlinear
Waves}, Wiley-Interscience.
\end{thebibliography}
\end{document}